

\hsize=14cm
\vsize=21cm
\parindent=0cm   \parskip=0pt
\pageno=1

\def\ind{\hskip 1cm\relax}

\ifnum\mag=\magstep1
\hoffset=-0.5cm   
\voffset=-0.5cm   
\fi

\pretolerance=500 \tolerance=1000  \brokenpenalty=5000

\catcode`\@=11

\font\eightrm=cmr8         \font\eighti=cmmi8
\font\eightsy=cmsy8        \font\eightbf=cmbx8
\font\eighttt=cmtt8        \font\eightit=cmti8
\font\eightsl=cmsl8        \font\sixrm=cmr6
\font\sixi=cmmi6           \font\sixsy=cmsy6
\font\sixbf=cmbx6


\font\tengoth=eufm10       \font\tenbboard=msbm10
\font\eightgoth=eufm8      \font\eightbboard=msbm8
\font\sevengoth=eufm7      \font\sevenbboard=msbm7
\font\sixgoth=eufm6        \font\fivegoth=eufm5

\skewchar\eighti='177 \skewchar\sixi='177
\skewchar\eightsy='60 \skewchar\sixsy='60


\newfam\gothfam           \newfam\bboardfam

\def\tenpoint{%
  \textfont0=\tenrm \scriptfont0=\sevenrm \scriptscriptfont0=\fiverm
  \def\rm{\fam\z@\tenrm}%
  \textfont1=\teni  \scriptfont1=\seveni  \scriptscriptfont1=\fivei
  \def\oldstyle{\fam\@ne\teni}\let\old=\oldstyle
  \textfont2=\tensy \scriptfont2=\sevensy \scriptscriptfont2=\fivesy
  \textfont\gothfam=\tengoth \scriptfont\gothfam=\sevengoth
  \scriptscriptfont\gothfam=\fivegoth
  \def\goth{\fam\gothfam\tengoth}%
  \textfont\bboardfam=\tenbboard \scriptfont\bboardfam=\sevenbboard
  \scriptscriptfont\bboardfam=\sevenbboard
  \def\bb{\fam\bboardfam\tenbboard}%
  \textfont\itfam=\tenit
  \def\it{\fam\itfam\tenit}%
  \textfont\slfam=\tensl
  \def\sl{\fam\slfam\tensl}%
  \textfont\bffam=\tenbf \scriptfont\bffam=\sevenbf
  \scriptscriptfont\bffam=\fivebf
  \def\bf{\fam\bffam\tenbf}%
  \textfont\ttfam=\tentt
  \def\tt{\fam\ttfam\tentt}%
  \abovedisplayskip=12pt plus 3pt minus 9pt
  \belowdisplayskip=\abovedisplayskip
  \abovedisplayshortskip=0pt plus 3pt
  \belowdisplayshortskip=4pt plus 3pt
  \smallskipamount=3pt plus 1pt minus 1pt
  \medskipamount=6pt plus 2pt minus 2pt
  \bigskipamount=12pt plus 4pt minus 4pt
  \normalbaselineskip=12pt
  \setbox\strutbox=\hbox{\vrule height8.5pt depth3.5pt width0pt}%
  \let\bigf@nt=\tenrm       \let\smallf@nt=\sevenrm
  \normalbaselines\rm}

\def\eightpoint{%
  \textfont0=\eightrm \scriptfont0=\sixrm \scriptscriptfont0=\fiverm
  \def\rm{\fam\z@\eightrm}%
  \textfont1=\eighti  \scriptfont1=\sixi  \scriptscriptfont1=\fivei
  \def\oldstyle{\fam\@ne\eighti}\let\old=\oldstyle
  \textfont2=\eightsy \scriptfont2=\sixsy \scriptscriptfont2=\fivesy
  \textfont\gothfam=\eightgoth \scriptfont\gothfam=\sixgoth
  \scriptscriptfont\gothfam=\fivegoth
  \def\goth{\fam\gothfam\eightgoth}%
  \textfont\bboardfam=\eightbboard \scriptfont\bboardfam=\sevenbboard
  \scriptscriptfont\bboardfam=\sevenbboard
  \def\bb{\fam\bboardfam}%
  \textfont\itfam=\eightit
  \def\it{\fam\itfam\eightit}%
  \textfont\slfam=\eightsl
  \def\sl{\fam\slfam\eightsl}%
  \textfont\bffam=\eightbf \scriptfont\bffam=\sixbf
  \scriptscriptfont\bffam=\fivebf
  \def\bf{\fam\bffam\eightbf}%
  \textfont\ttfam=\eighttt
  \def\tt{\fam\ttfam\eighttt}%
  \abovedisplayskip=9pt plus 3pt minus 9pt
  \belowdisplayskip=\abovedisplayskip
  \abovedisplayshortskip=0pt plus 3pt
  \belowdisplayshortskip=3pt plus 3pt
  \smallskipamount=2pt plus 1pt minus 1pt
  \medskipamount=4pt plus 2pt minus 1pt
  \bigskipamount=9pt plus 3pt minus 3pt
  \normalbaselineskip=9pt
  \setbox\strutbox=\hbox{\vrule height7pt depth2pt width0pt}%
  \let\bigf@nt=\eightrm     \let\smallf@nt=\sixrm
  \normalbaselines\rm}

\def\pc#1{\bigf@nt#1\smallf@nt}         \def\pd#1 {{\pc#1} }
\catcode`\;=\active
\def;{\relax\ifhmode\ifdim\lastskip>\z@\unskip\fi
\kern\fontdimen2  -1.2 \fontdimen3 \string;}

\catcode`\:=\active
\def:{\relax\ifhmode\ifdim\lastskip>\z@\unskip\fi\penalty\@M\ \fi\string:}

\catcode`\!=\active
\def!{\relax\ifhmode\ifdim\lastskip>\z@
\unskip\fi\kern\fontdimen2  -1.1 \fontdimen3 \string!}

\catcode`\?=\active
\def?{\relax\ifhmode\ifdim\lastskip>\z@
\unskip\fi\kern\fontdimen2  -1.1 \fontdimen3 \string?}

\def\^#1{\if#1i{\accent"5E\i}\else{\accent"5E #1}\fi}
\def\"#1{\if#1i{\accent"7F\i}\else{\accent"7F #1}\fi}

\frenchspacing

\newtoks\auteurcourant      \auteurcourant={\hfil}
\newtoks\titrecourant       \titrecourant={\hfil}

\newtoks\hautpagetitre      \hautpagetitre={\hfil}
\newtoks\baspagetitre       \baspagetitre={\hfil}

\newtoks\hautpagegauche
\hautpagegauche={\eightpoint\rlap{\folio}\hfil\the\auteurcourant\hfil}
\newtoks\hautpagedroite
\hautpagedroite={\eightpoint\hfil\the\titrecourant\hfil\llap{\folio}}

\newtoks\baspagegauche      \baspagegauche={\hfil}
\newtoks\baspagedroite      \baspagedroite={\hfil}

\newif\ifpagetitre          \pagetitretrue

\headline={\ifpagetitre\the\hautpagetitre
\else\ifodd\pageno\the\hautpagedroite\else\the\hautpagegauche\fi\fi}

\footline={\ifpagetitre\the\baspagetitre\else
\ifodd\pageno\the\baspagedroite\else\the\baspagegauche\fi\fi
\global\pagetitrefalse}
\def\raggedbottom{\topskip 10pt plus 36pt\r@ggedbottomtrue}
\def\pointir{\unskip . --- \ignorespaces}
\def\Bigbreak{\vskip-\lastskip\bigbreak}
\def\Medbreak{\vskip-\lastskip\medbreak}
\def\ctexte#1\endctexte{%
  \hbox{$\vcenter{\halign{\hfill##\hfill\crcr#1\crcr}}$}}
\long\def\ctitre#1\endctitre{%
    \ifdim\lastskip<24pt\vskip-\lastskip\bigbreak\bigbreak\fi
  		\vbox{\parindent=0pt\leftskip=0pt plus 1fill
          \rightskip=\leftskip
          \parfillskip=0pt\bf#1\par}
    \bigskip\nobreak}

\long\def\section#1\endsection{%
\vskip 0pt plus 3\normalbaselineskip
\penalty-250
\vskip 0pt plus -3\normalbaselineskip
\Bigbreak
\message{[section \string: #1]}{\bf#1\unskip}\pointir}

\long\def\sectiona#1\endsection{%
\vskip 0pt plus 3\normalbaselineskip
\penalty-250
\vskip 0pt plus -3\normalbaselineskip
\Bigbreak
\message{[sectiona \string: #1]}%
{\bf#1}\medskip\nobreak}

\let\+=\tabalign

\def\signature#1\endsignature{\vskip 15mm minus 5mm\rightline{\vtop{#1}}}

\long\def\th#1 #2\enonce#3\endth{%
   \Medbreak
   {\pc#1} {#2\unskip}\pointir{\it #3}\medskip}

\long\def\tha#1 #2\enonce#3\endth{%
   \Medbreak
   {\pc#1} {#2\unskip}\par\nobreak{\it #3}\medskip}

\def\decale#1{\smallbreak\hskip 28pt\llap{#1}\kern 5pt}
\def\decaledecale#1{\smallbreak\hskip 34pt\llap{#1}\kern 5pt}
\def\puce{\smallbreak\hskip 6pt{$\scriptstyle\bullet$}\kern 5pt}
\mathcode`A="7041 \mathcode`B="7042 \mathcode`C="7043 \mathcode`D="7044
\mathcode`E="7045 \mathcode`F="7046 \mathcode`G="7047 \mathcode`H="7048
\mathcode`I="7049 \mathcode`J="704A \mathcode`K="704B \mathcode`L="704C
\mathcode`M="704D \mathcode`N="704E \mathcode`O="704F \mathcode`P="7050
\mathcode`Q="7051 \mathcode`R="7052 \mathcode`S="7053 \mathcode`T="7054
\mathcode`U="7055 \mathcode`V="7056 \mathcode`W="7057 \mathcode`X="7058
\mathcode`Y="7059 \mathcode`Z="705A

\def\spacedmath#1{\def\packedmath##1${\bgroup\mathsurround=0pt ##1\egroup$}%
\mathsurround#1 \everymath={\packedmath}\everydisplay={\mathsurround=0pt }}

\def\nospacedmath{\mathsurround=0pt \everymath={}\everydisplay={} }

\def\up#1{\raise 1ex\hbox{\smallf@nt#1}}

\def\qed{\raise -2pt\hbox{\vrule\vbox to 10pt{\hrule width 4pt
                 \vfill\hrule}\vrule}}

\def\virg{\raise .4ex\hbox{,}}   


\def\build#1_#2^#3{\mathrel{
\mathop{\kern 0pt#1}\limits_{#2}^{#3}}}
\catcode`\@=12

\showboxbreadth=-1  \showboxdepth=-1

\magnification=1200

\tenpoint

\spacedmath{2pt}
\baselineskip 14pt
\font\ninesl=cmsl9
\font\sevensl=cmsl10 scaled 700
\vskip 2cm
 \ctitre
Sur la cohomologie de certains espaces de modules de fibr\'es
vectoriels
 \endctitre
\centerline{Arnaud {\pc BEAUVILLE}\footnote{(*)}{\sevenrm Avec le support
partiel
du projet europ\'een Science ``Geometry of Algebraic Varieties", Contrat
n\up{o} SCI-0398-C (A). adresse \'electronique: arnaud@matups.matups.fr}}
 \vskip 1cm
\hskip 8cm\leftline{\ninesl D\'edi\'e \`a M.S. Narasimhan et C.S. Seshadri,}

\hskip 8cm\leftline {\ninesl pour leur 60\up{\sevensl\`eme} anniversaire.}
\vskip 1cm
 \ind Soit $X$ une surface de Riemann
compacte. Fixons des entiers $r$ et $d$ {\it premiers entre eux}, avec $r\ge
1$, et notons ${\cal M}$ l'espace des modules ${\cal U}_X(r,d)$ des fibr\'es
vectoriels stables sur $X$, de rang $r$ et de degr\'e $d$. C'est une
vari\'et\'e
projective et lisse, et il existe un {\it fibr\'e de Poincar\'e} ${\cal E}$ sur
$X \times {\cal M}$ ; cela signifie que pour tout point $e$ de ${\cal M}$,
correspondant \`a un fibr\'e $E$ sur $X$, la restriction de ${\cal E}$ \`a $X
\times \{ e\}$ est isomorphe \`a $E$.

\ind Notons $p$ , $q$ les projections de $X \times {\cal M}$ sur $X$ et
${\cal M}$ respectivement. Soit $m$ un entier $\le r$ ; la classe de Chern
$c_m({\cal E})$ admet une d\'ecomposition de K\"unneth
$$c_m({\cal E})=\sum_i p^*\xi_i \cdot q^*\mu_i\ ,\quad {\rm avec}\quad \xi_i
\in H^*(X,{\bf
Z})\  ,\ \mu_i \in H^*({\cal M},{\bf Z})\ ,\quad {\rm deg}(\xi_i) +{\rm
deg}(\mu_i) = 2m .$$

Nous dirons que les classes $\mu_i$ sont les composantes de K\"unneth de
$c_m({\cal E})$. Un des r\'esultats essentiels de [A--B] est la d\'etermination
d'un ensemble de g\'en\'erateurs de l'alg\`ebre de cohomologie $H^*({\cal
M},{\bf
Z})$ ; il a la cons\'equence suivante :

\th TH\'EOR\`EME
\enonce L'alg\`ebre de cohomologie $H^*({\cal M},{\bf Q})$
est engendr\'ee par les composantes de K\"unneth des classes de Chern de
${\cal E}$ .
\endth
\ind  Le but de cette note est de montrer comment la m\'ethode de
la diagonale utilis\'ee dans [E--S] fournit une d\'emonstration tr\`es simple
de ce th\'eor\`eme. Celui-ci r\'esulte de l'\'enonc\'e un peu plus g\'en\'eral
que voici:

\th PROPOSITION
\enonce Soit $X$ une vari\'et\'e complexe projective et
lisse, et ${\cal M}$ un espace de modules de faisceaux stables sur $X$ (par
rapport \`a une polarisation fix\'ee, cf. {\rm [M]}). On fait les hypoth\`eses
suivantes :

\decale{(i)} La vari\'et\'e ${\cal M}$ est projective et lisse.
\decale{(ii)} Il existe un faisceau de Poincar\'e ${\cal E}$ sur ${\cal M}$.
\decale{(iii)} Pour $E$ , $F$ dans ${\cal M}$ , on a ${\rm Ext}^i(E,F) = 0$
pour $i\ge
2$ .

Alors l'alg\`ebre de cohomologie $H^*({\cal M},{\bf Q})$
est engendr\'ee par les composantes de K\"unneth des classes de Chern de
${\cal E}$ .
 \endth
\ind La d\'emonstration suit de pr\`es celle du th. 1 de [E--S]. Rappelons-en
l'id\'ee
fondamentale: soit $\delta$ la classe de cohomologie de la
diagonale dans $H^*({\cal M}\times {\cal M},{\bf Q})$; notons $p$ et $q$ les
deux
projections de  ${\cal M}\times {\cal M}$ sur ${\cal M}$. Soit $\displaystyle
\delta =
\sum_i p^*\mu_i \cdot q^*\nu_i$ la d\'ecomposition de K\"unneth de $\delta$;
alors {\it
l'espace $H^*({\cal M},{\bf Q})$ est engendr\'e par les $\nu_i$.} En effet,
pour toute classe
$\lambda$ dans $H^*({\cal M},{\bf Q})$ , on a
$$\lambda = q_*(\delta\cdot p^*\lambda) = \sum {\rm deg}(\lambda\cdot
\mu_i)\nu_i\quad ,$$
d'o\`u notre assertion. Il s'agit donc d'exprimer la classe $\delta$ en
fonction des classes de
Chern du fibr\'e universel.

\ind Notons
$p_1,p_2$ les deux projections de $C\times {\cal M}\times {\cal M}$ sur
$C\times {\cal M}$ ,
$\pi$ la projection sur ${\cal M}\times {\cal M}$ ; d\'esignons par ${\cal H}$
le faisceau ${\cal
H}om(p_1^*{\cal E},p_2^*{\cal E})$ . Vu l'hypoth\`ese (iii), l'hypercohomologie
$R{\pi}_*{\cal H}$
est repr\'esent\'ee dans la cat\'egorie d\'eriv\'ee par un complexe de fibr\'es
$K^{\bullet}$ , nul
en degr\'e diff\'erent de $0$ et $1$. Autrement dit, il existe un morphisme de
fibr\'es
$u:K^0\longrightarrow K^1$ tel qu'on ait, pour tout point $x=(E,F)$ de ${\cal
M}$, une suite exacte
\nospacedmath
$$ 0\rightarrow {\rm Hom}(E,F)\longrightarrow K^0(x) \buildrel{
u(x)}\over{\hbox to 1cm{\rightarrowfill}} K^1(x) \longrightarrow
{\rm Ext}^1(E,F)\rightarrow 0 \quad .$$
\spacedmath{2pt}
 \ind Comme l'espace ${\rm Hom}(E,F)$ est non nul si et seulement si $E$ et $F$
sont isomorphes, on voit que la diagonale $\Delta$ de ${\cal M}\times {\cal
M}$ co\"incide ensemblistement avec le lieu de
d\'eg\'en\'erescence $D$ de $u$ (d\'efini par l'annulation des mineurs de rang
maximal de $u$). On peut prouver comme dans [E--S] l'\'egalit\'e sch\'ematique,
mais cela n'est
pas n\'ecessaire pour d\'emontrer la proposition.

\ind Soit $E$ un \'el\'ement de ${\cal M}$ . On a
$$
{\rm rg}(K^0)-{\rm rg}(K^1) = {\rm dim\ Hom}(E,F)\,-\,{\rm dim\ Ext}^1(E,F)
$$
quel que soit le point $(E,F)$ de ${\cal M}\times {\cal M}$ .
Puisque ${\rm Ext}^2(E,E)=0$ , la dimension $m$ de ${\cal M}$ est
\'egale \`a ${\rm dim\ Ext}^1(E,E)$ ; ainsi la sous-vari\'et\'e
d\'eterminantale $D$ de ${\cal M}\times {\cal M}$ a la codimension attendue
${\rm rg}(K^1)-{\rm rg}(K^0)+1$ . Sa classe de cohomologie $\delta '\in
H^m({\cal M}\times {\cal M},{\bf Z})$  est alors donn\'ee par la formule de
Porteous $$ \delta ' = c_m(K^1-K^0) = c_m(-{\pi}_!{\cal H})\ , $$
o\`u ${\pi}_!$ d\'esigne le foncteur image directe en K-th\'eorie. Cette classe
\'etant multiple
de la classe $\delta$ de la diagonale, on conclut avec le lemme suivant:

\th LEMME
\enonce Soit ${\cal A}$ la sous-${\bf Q}$-alg\`ebre de $H^*({\cal M},{\bf Q})$
engendr\'ee par les composantes de K\"unneth des classes de Chern de
${\cal E}$ , et soient $p$ et $q$ les deux projections de  ${\cal M}\times
{\cal M}$ sur ${\cal M}$. Les classes de Chern de ${\pi}_!{\cal H}$ sont de
la forme $\sum p^*\mu_i \cdot q^*\nu_i$ , avec $\mu_i\, ,\,\nu_i \in
{\cal A} . $
\endth \ind Notons $r$ la projection de $C\times {\cal M}\times {\cal
M}$ sur $C$ . Tout polyn\^ome en les classes de Chern de $p_1^*{\cal E}$ et de
$p_2^*{\cal E}$ est une somme de produits de la forme $r^*\gamma \cdot
{\pi}^*p^*\mu \cdot {\pi}^*q^*\nu$ , o\`u $\mu$ et $\nu$
appartiennent \`a ${\cal A}$ . Le lemme r\'esulte alors de la formule de
Riemann-Roch
$$ ch({\pi}_!{\cal H}) = {\pi}_*\bigl(r^*{\rm Todd(C)\,ch}({\cal H})\bigr)
\quad .
\quad\vrule height 4pt depth 0pt width 4pt$$

\bigskip
{\it Remarque}.-- La condition (iii) de la proposition est \'evidemment
tr\`es contraignante. Donnons deux exemples:
\decale{a)} $X$ est une surface rationnelle ou r\'egl\'ee, et la polarisation
$H$ v\'erifie $H\cdot K_X<0$. L'argument de [M, cor. 6.7.3] montre que la
condition (iii) est satisfaite. Si de plus les coefficients $a_i$ du
polyn\^ome de Hilbert des \'el\'ements de ${\cal M}$, \'ecrit sous la forme
$\displaystyle \chi (E\otimes H^m) = \sum_{i=0}^2 a_i{m+i\choose i}$ , sont
premiers entre eux, les conditions (i) \`a (iii) sont satisfaites [M, \S 6].

\ind Dans le cas d'une surface {\it rationnelle}, on obtient mieux. Pour
toute vari\'et\'e $T$, d\'esi\nobreak gnons par $CH^*(T)$ l'anneau de Chow de
$T$;
gr\^ace \`a l'isomorphisme $CH^*(X\times {\cal M}) \cong$\penalty -10000 $
CH^*(X)\otimes
CH^*({\cal M})$, on peut remplacer dans la d\'emonstration de la proposition
l'anneau de cohomologie par l'anneau de Chow. On en d\'eduit que {\it la
cohomologie rationnelle de ${\cal M}$ est alg\'ebrique,} c'est-\`a-dire que
l'application "classe de cycles" $CH^*({\cal M})\otimes {\bf
Q}\longrightarrow H^*({\cal M},{\bf Q})$ est un isomorphisme d'anneaux. Dans
le cas $X={\bf P}^2$, Ellingsrud et Str\o mme obtiennent le m\^eme r\'esultat
{\it sur} ${\bf Z}$, plus le fait que ces groupes sont sans torsion, gr\^ace
\`a
l'outil suppl\'ementaire de la suite spectrale de Beilinson.

\decale{b)} $X$ est une vari\'et\'e de Fano de dimension $3$. Soit
$S$ une surface lisse appartenant au syst\`eme lin\'eaire
$|-K_X|$ (de sorte que $S$ est une surface K3). Lorsqu'elle est satisfaite, la
condition (iii) a des cons\'equences remarquables [T]: elle entra\^ine que
l'application de restriction $E \longmapsto E_{|S}$ d\'efinit un isomorphisme
de ${\cal M}$ sur une {\it sous-vari\'et\'e lagrangienne} d'un espace de
modules ${\cal M}_S$ de fibr\'es sur $S$ (muni de sa structure
symplectique canonique). Il me semble int\'eressant de mettre en
\'evidence des espaces de modules de fibr\'es sur une vari\'et\'e de Fano (et
d\'ej\`a sur ${\bf P}^3$) poss\'edant la propri\'et\'e (iii).

\vskip 2cm
\centerline{BIBLIOGRAPHIE}
\bigskip
\hangindent 3cm
\hangafter 1
[A--B]\quad M. {\pc ATIYAH} et R. {\pc BOTT} : {\it The Yang-Mills equations
over
Riemann surfaces.} Phil. Trans. R. Soc. London A 308, 523-615 (1982).
\medskip
\hangindent 3cm
\hangafter 1
[E--S]\quad G. {\pc ELLINGSRUD} et S.A. {\pc STR\O MME} : {\it On generators
for the Chow
ring of fine moduli spaces on ${\bf P}^2$}. Preprint (1989).
\medskip
\hangindent 3cm
\hangafter 1
[M]\quad M. {\pc MARUYAMA} : {\it Moduli of
stable sheaves, II.} J. Math. Kyoto Univ. 18 (1978), 557-614.
\medskip
\hangindent 3cm
\hangafter 1
[T]\quad A. N. {\pc TYURIN} : {\it The moduli space of vector bundles on
threefolds,
surfaces and curves I.} Preprint (1990).

\bigskip
\signature
\+ Arnaud {\pc BEAUVILLE}\cr
\+ Universit\'e Paris-Sud\cr
\+ Math\'ematiques -- B\^at. 425\cr
\+ 91\hskip 3pt 405 Orsay Cedex, France\cr
\endsignature
\bye